\documentclass[12pt]{amsart}

\usepackage{epsfig}
\usepackage{graphics}
\usepackage{latexsym}
\usepackage{psfrag}
\usepackage{amssymb}
\usepackage{amsmath}
\usepackage{amsthm} 

\renewcommand{\setminus}{{\smallsetminus}}

\newcommand{\Khat}{\hat{K}}
\newcommand{\tauhat}{\hat{\tau}}
\newcommand{\ahat}{\hat{a}}

\newcommand{\bdy}{{\partial}} % Boundary
\newcommand{\map}[1]{\overset{#1}{\longrightarrow}}

\theoremstyle{plain}
\newtheorem{theorem}{Theorem}[section]
\newtheorem{corollary}[theorem]{Corollary}
\newtheorem{lemma}[theorem]{Lemma}
\newtheorem{claim}[theorem]{Claim}

\theoremstyle{definition}
\newtheorem*{define}{Definition}
\newtheorem*{remark}{Remark}

\newtheorem*{notation}{Notation}

\newsavebox{\savepar}

\begin{document}

\title{Involutions of Knots That Fix Unknotting Tunnels} 
\address{Mathematics Department, Stanford University}
\author{David Futer}
\date{\today}

\begin{abstract}
Let $K$ be a knot that has an unknotting tunnel
$\tau$. We prove that $K$ admits a strong involution that fixes $\tau$
pointwise if and only if $K$ is a two-bridge knot and $\tau$ its
upper or lower tunnel.
\end{abstract}

\maketitle

\section{Introduction}

Let $L$ be a link of one or two components in $S^3$. An {\em
unknotting tunnel} for $L$ is a properly embedded arc $\tau$, with
$L \cap \tau = \bdy \tau$, such that the complement of a regular
neighborhood of $L \cup \tau$ is a genus $2$ handlebody. As described
in \cite{bm}, an unknotting tunnel induces a {\em strong inversion} of
the link complement -- that is, an involution of $S^3$ that sends each
component of $L$ to itself with reversed orientation.

When $L$ is a two-component link, it is known that this involution can
be chosen to fix $\tau$ pointwise. As a result, an argument of Adams
in \cite{adams} gives us some geometric information about the tunnel:
when the complement of $L$ is hyperbolic, $\tau$ is isotopic to a
geodesic in the geometric structure. \cite{adams} conjectures that the
same is true for knots: that any unknotting tunnel for a hyperbolic
knot is isotopic to a geodesic. As for links, this conjecture would
follow easily if it were known that $\tau$ is fixed pointwise by some
strong inversion of the knot.

The main result of this paper is that this happens only in the well-known
special case of $2$-bridge knots. (See Section \ref{2bridge} for
background on on $2$-bridge knots and their unknotting tunnels.)

\begin{theorem}\label{main}
Let $K \subset S^3$ be a knot with an unknotting tunnel $\tau$. Then
$\tau$ is fixed pointwise by a strong inversion of $K$ if and only if
$K$ is a two-bridge knot and $\tau$ is its upper or lower tunnel.
\end{theorem}

This theorem materialized amid the wonderful hospitality of the Isaac
Newton Institute for Mathematical Sciences at Cambridge University. I
have benefited greatly from conversations with Ian Agol, Yoav Moriah,
and Makoto Sakuma. Steve Kerckhoff and Saul Schleimer both devoted a
great deal of time to hearing and vetting many versions of the
proof-in-progress, and deserve my sincerest gratitude.

\section{Tunnels and Involutions}

\begin{notation}
From now on, $K$ will denote a knot in $S^3$ that has an unknotting
tunnel $\tau$. Thus $K \cup \tau$ is realized as a graph with two
vertices and three edges $\tau$, $K_1$, and $K_2$, where $K = K_1 \cup
K_2$. When identifying $\tau$ as a particular tunnel $\tau_0$ of the
knot $K$, we mean that $K \cup \tau$ is equivalent to $K \cup \tau_0$
via an isotopy of $S^3$ that preserves $K$ setwise. Note that this
notion of equivalence is stronger than isotopy of tunnels in the knot
exterior, because the endpoints of $\tau$ are not allowed to pass
through each other.
\end{notation}

\begin{define}
For any graph $\Gamma \subset S^3$, let $N(\Gamma)$ be an open regular
neighborhood of $\Gamma$. We call $E(\Gamma) = S^3
\setminus N(\Gamma)$ the {\em exterior} of $\Gamma$.
\end{define}

An unknotting tunnel induces a genus $2$ Heegaard splitting of $S^3$
into handlebodies $V_1 = \overline{N(K \cup \tau)}$ and $V_2 = E(K
\cup \tau)$. Let $\Sigma$ be the Heegaard surface: $\Sigma = \bdy 
N(K \cup \tau)$. A genus $2$ handlebody admits a {\em hyper-elliptic
involution} that preserves the isotopy class of every simple closed
curve on its boundary; it is unique up to isotopy. Thus the
hyper-elliptic involutions of $V_1$ and $V_2$ can be joined over
$\Sigma$ to an orientation-preserving involution $\varphi_\tau$ on all
of $S^3$. This involution sends $K$ to itself with reversed
orientation; i.e. is a {\em strong inversion} of $K$ (see \cite{bm}).

The action of $\varphi_\tau$ preserves the meridians of $K_1$ and
$K_2$ on $\Sigma$ while reversing the orientation on $K$ -- so it must
switch the two vertices of $K \cup \tau$ and thus reverse the
orientation on $\tau$. It is conceivable, however, that some other
involution $\psi$ might fix $\tau$ pointwise while switching $K_1$
with $K_2$. Theorem \ref{main} says that this only happens for
$2$-bridge knots.

\section{Two-Bridge Knots}\label{2bridge}

\begin{define} A {\em rational tangle} is a pair $(B, t)$, where $B$ is 
a 3-ball and $t$ consists of two disjoint, properly embedded arcs
$\gamma_1$ and $\gamma_2$. We further require that the $\gamma_i$ are
both isotopic to $\bdy B$ via disjoint disks $D_i$.
\end{define}

\begin{define}
A knot $K \subset S^3$ is called a {\em two-bridge knot} if some
sphere $S$ splits $(S^3, K)$ into two rational tangles. $S$ is then
called a {\em bridge sphere} for $K$.
\end{define}

A $2$-bridge knot has a {\em 4-plat projection} with two maxima at the
top and two minima at the bottom \cite{bz}. Any horizontal plane,
together with the point at infinity, then serves as a bridge sphere
for $K$. The $4$-plat projection also reveals two unknotting tunnels
for $K$: an {\em upper tunnel} $\tau_1$ connecting the two maxima, and
a {\em lower tunnel} $\tau_2$ connecting the two minima. 

Two-bridge knots also have four {\em dual tunnels}, which may be
isotopic to the upper or lower tunnels in special cases. Morimoto and
Sakuma classified the six total tunnels up to homeomorphism and
isotopy of the knot exterior \cite{ms}. Kobayashi showed that any
unknotting tunnel of a $2$-bridge knot is exterior-isotopic to one of
these six \cite{kob}.

The following result is an immediate consequence of \cite[Section 3]{bm}.

\begin{lemma}\label{known-case}
Let $K$ be a two-bridge knot, and $\tau_1$ and $\tau_2$ be its upper
and lower tunnels. Let $\varphi_{\tau_i}$ the involution of $S^3$ that
comes from the handlebody decomposition induced by $\tau_i$. Then
$\varphi_{\tau_1}$ can be chosen to fix $\tau_2$ pointwise, and
$\varphi_{\tau_2}$ can be chosen to fix $\tau_1$ pointwise.
\end{lemma}

\begin{figure}
\psfrag{tau1}{$\tau_1$}
\psfrag{tau2}{$\tau_2$}
\psfrag{ft1}{$\varphi_{\tau_1}$}
\psfrag{ft2}{$\varphi_{\tau_2}$}
\begin{center}
\epsfbox{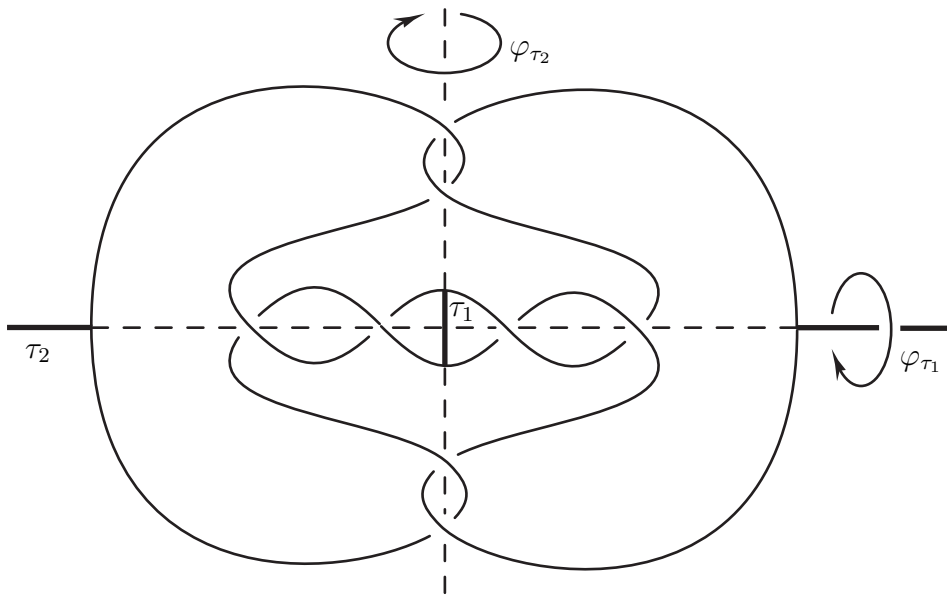}
\caption{Upper and lower tunnels, and the corresponding involutions, in 
the tri-symmetric projection.}
\label{tri-sym}
\end{center}
\end{figure}

\begin{proof}
As described in \cite{bz}, there is an isotopy of $K$ from the
$4$-plat projection to a {\em tri-symmetric projection}, carrying
$\tau_1$ to part of the vertical axis of symmetry and $\tau_2$ to part
of the horizontal axis. (See Figure \ref{tri-sym}.) Abusing notation
slightly, we continue to call these isotoped tunnels $\tau_1$ and
$\tau_2$. Now, the involution $\varphi_{\tau_1}$ is evident in the
figure as a $180^\circ$ rotation about the horizontal axis containing
$\tau_2$, and thus fixes $\tau_2$ pointwise. Similarly,
$\varphi_{\tau_2}$ is a rotation about the vertical axis containing
$\tau_1$, fixing $\tau_1$ pointwise.
\end{proof}

For the purpose of recognizing two-bridge knots based on involutions,
our main tool is a theorem of Scharlemann and Thompson about
planar graphs.

\begin{define}
A finite graph $\Gamma \subset S^3$ is called {\em planar} if it lies
on an embedded $2$-sphere in $S^3$. $\Gamma$ is called {\em abstractly
planar} if it is homeomorphic to a planar graph.
\end{define}

\begin{theorem}\cite{st}\label{graph-thm}
A finite graph $\Gamma \subset S^3$ is planar if and only if
\begin{enumerate}
\item $\Gamma$ is abstractly planar,
\item every proper subgraph of $\Gamma$ is planar, and
\item $E(\Gamma)$ is a handlebody.
\end{enumerate}
\end{theorem}

\begin{remark}
For our purposes, we will only need the special case of this theorem
when $\Gamma$ consists of one vertex and two edges. This was first
proved in an unpublished preprint of Hempel and Roeling.
\end{remark}

\begin{corollary}\label{recognition}
Let $K$ be a knot with unknotting tunnel $\tau$, where $\tau$ splits
$K$ into edges $K_1$ and $K_2$. If both $K_1 \cup \tau$ and $K_2 \cup
\tau$ are unknots, then $K$ is a two-bridge knot with splitting sphere 
$\bdy N(\tau)$. Furthermore, $\tau$ is an upper or lower tunnel.
\end{corollary}

\begin{proof}
Let $f:S^3 \to S^3$ be a map that contracts $\overline{N(\tau)}$ to a
point and is the identity outside a small regular neighborhood of
$\overline{N(\tau)}$. Then $\Gamma = f(K \cup \tau)$ is a graph with
one vertex and two edges, which is clearly abstractly planar. Each
proper subgraph $\Gamma_i \subset \Gamma$ is the image of $K_i \cup
\tau$, and is thus an unknotted circle. Also, $E(\Gamma)$ is a
handlebody since $E(K \cup \tau)$ is a handlebody. Thus, by Theorem
\ref{graph-thm}, $\Gamma$ is planar.

\begin{figure}
\psfrag{K1}{$K_1$}
\psfrag{K2}{$K_2$}
\psfrag{E1}{$E_1$}
\psfrag{E2}{$E_2$}
\psfrag{G1}{$\Gamma_1$}
\psfrag{G2}{$\Gamma_2$}
\psfrag{D1}{$D_1$}
\psfrag{D2}{$D_2$}
\psfrag{f}{$f$}
\psfrag{bnt}{$\bdy N(\tau)$}
\begin{center}
\epsfbox{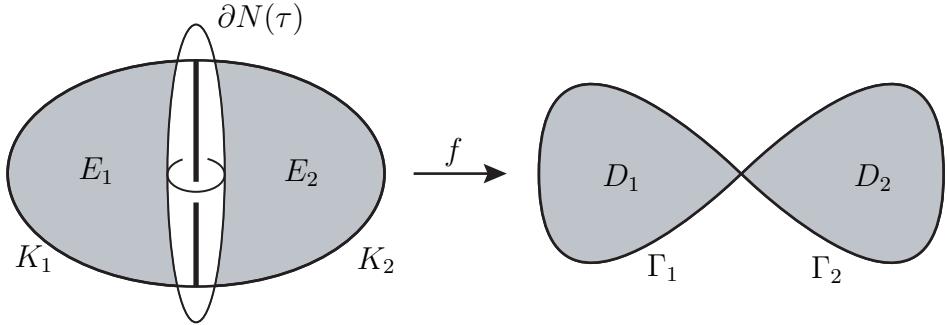}
\caption{The graphs and disks of Corollary \ref{recognition}.}
\label{planar-fig}
\end{center}
\end{figure}

Since $\Gamma$ is planar, its subgraphs $\Gamma_1$ and $\Gamma_2$
bound disjoint disks $D_1$ and $D_2$ in $S^3$. (See Figure
\ref{planar-fig}.) These disks pull back via $f$ to disjoint disks
$E_1, E_2 \subset E(\tau)$. Each $E_i$ thus provides an isotopy of
$K_i \cap E(\tau)$ to $\bdy N(\tau)$, making $(E(\tau), K \cap
E(\tau))$ a rational tangle. Meanwhile, $N(\tau)$ intersects $K$ in
two short arcs that are clearly boundary-parallel, so that is a
rational tangle too. Therefore, $K$ is a $2$-bridge knot with
splitting sphere $\bdy N(\tau)$.

It is known (from \cite{schubert}, for example) that a splitting of a
$2$-bridge knot into rational tangles is unique up to isotopy. Thus an
ambient isotopy of $S^3$ carries $K$ to a $4$-plat projection and
$N(\tau)$ to a half-space containing just the two maxima (or
minima). Thus $\tau$ is an upper or lower tunnel for $K$.
\end{proof}

\section{Cyclic Groups Acting on Handlebodies}

A key ingredient in the proof of Theorem \ref{main} is the equivariant
loop theorem of Meeks and Yau.

\begin{theorem}\cite{my}\label{elt}
Let $M$ be a compact, connected $3$-manifold with connected boundary,
and $G$ be a finite group acting smoothly on $M$. Let $H$ be the
kernel of the homomorphism of fundamental groups induced by the
inclusion $\bdy M \hookrightarrow M$. Then there is a collection
$\Delta = \{D_1, \ldots, D_n\}$ of disjoint, properly embedded,
essential disks in $M$ such that

\begin{enumerate}
\item $\bdy D_i \subset \bdy M$ for all $i$,
\item $H$ is the normal subgroup of $\pi_1(M)$ generated by the $\bdy D_i$, and
\item $\Delta$ is invariant under the action of $G$.
\end{enumerate}
\end{theorem}

\begin{notation}
Let $M$ and $\Delta$ be as above, and denote by $S$ be the union of
the disks in $\Delta$. Let $M \vert \Delta = M \setminus N(S)$, where
the regular neighborhood is chosen to be invariant under the action of
$G$.
\end{notation}

To adapt this theorem for our purposes, we assume that $M$ is a
handlebody and analyze the pieces of $M \vert \Delta$.

\begin{lemma}\label{ball-pieces}
Let $V$ be a handlebody, and let $G$ be a finite group acting smoothly
on $V$. Let $\Delta$ be the collection of disks guaranteed by Theorem
\ref{elt}. Then each component of $V \vert \Delta$ is a ball.
\end{lemma}

\begin{proof}
First, suppose that $V$ is cut along just one disk $D$. It is
well-known that the result is one or two handlebodies. Here is an
outline of the proof, suggested by Saul Schleimer. Let $V'$ be one
component of the resulting manifold. Then $\pi_1(V')$ injects into
$\pi_1(V)$, because if a loop $\gamma \subset V'$ bounds a disk $E
\subset V$, one can do disk swaps with $D$ to find a disk that $\gamma$ 
bounds in $V'$ also. Thus $\pi_1(V')$ is free (as a subgroup of the
free group $\pi_1(V)$), and so $V'$ is a handlebody. Applying this
argument repeatedly, we see that every component of $V \vert \Delta$
is a handlebody.

Let $X$ be one component of $V \vert \Delta$. $\bdy X$ consists of a
subset of $\bdy V$ along with some number of distinguished disks $E_1,
\ldots, E_k$ that come from removing $N(S)$. We already know that $X$
is a handlebody, and we prove that $X$ is a ball by showing that it
contains no essential disks.

Suppose that $D \subset X$ is a disk with $\gamma = \bdy D \subset
\bdy X$. Clearly, $D$ can be isotoped so that $\gamma$ is disjoint from
the distinguished disks $E_j$. Thus $\gamma$ is a simple closed curve
in $\bdy V$ that bounds a disk $D \subset V$. By Theorem \ref{elt},
some loop freely homotopic to $\gamma$ is generated by conjugates of
the $\bdy D_i$ in $\pi_1(\bdy V)$. Passing to homology, we see that
the class $[\gamma]$ is a linear combination of the $[\bdy D_i]$ in
$H_1(\bdy V)$.

Let us express $\bdy V$ as a union of open subsets $A$ and $B$, where
$B = \bdy V \setminus \bdy X$ and $A$ is an open regular
neighborhood of $\bdy X$ in $\bdy V$. Thus $A \cap B$ is a disjoint
union of open regular neighborhoods of the $\bdy E_j$. $\gamma$ lies in 
$A$ but is homologous to some cycle $c$ in $B$, because

\vspace{-0.2in}

$$[\gamma] = \sum_{i=1}^n a_i [\bdy D_i] \in H_1(\bdy V) 
\mbox{, and } \bdy D_i \subset B \mbox{ for all } i.$$

Now, the Mayer-Vietoris sequence gives us

\vspace{-0.2in}
$$\ldots \to H_1(A \cap B) \map{i_\ast} H_1(A) \oplus H_1(B) \map{j_\ast} 
H_1(\bdy V) \to \ldots$$

\noindent induced by the chain maps $i(x) = (x, x)$ and $j(x,y) = x-y$.
Since $[\gamma] - [c] = 0 \in H_1(\bdy V)$, $([\gamma], [c]) \in
\ker(j_\ast)$. But since the sequence is exact, $([\gamma], [c]) \in
\mathrm{Im}(i_\ast)$. Thus $\gamma$ is homologous in $H_1(A)$ to a
cycle that lies in $A \cap B$. But $H_1(A \cap B)$ is generated by the
$[\bdy E_j]$, so $\gamma$ must be homologically trivial in $\bdy X$.

Since $X$ is a handlebody without any homologically essential disks,
it must be a ball.
\end{proof}

The following application of Meeks and Yau's theorem to a form of the
Smith Conjecture for handlebodies was suggested by Ian Agol and Saul
Schleimer.

\begin{theorem}\label{smith-hand}
Let $V$ be a handlebody, and let $g:V \to V$ be an orientation-preserving, 
periodic diffeomorphism. Then

\begin{enumerate}
\item The fixed-point set of $g$ is either empty or boundary-parallel. 
\item $W = V/(g)$ is also a handlebody, in which the image of the 
fixed-point set is again empty or boundary-parallel.
\end{enumerate}
\end{theorem}

\begin{proof}
Let $\Delta$ be the collection of disks in $V$ given by Theorem
\ref{elt}. By Lemma \ref{ball-pieces}, $V \vert \Delta$ is a disjoint 
union of balls. The fixed-point set $a$ may be empty, or it may have
one or more components. If $a$ is non-empty, we will prove that it is
boundary-parallel by first considering its intersection with the
individual balls and then seeing how the pieces join up.

For each component $B$ of $V \vert \Delta$, $\bdy B$ consists of a
subset of $\bdy V$, together with some number of distinguished disks
$E_1, \ldots, E_k$ that come from removing $N(S)$. Since $V \vert
\Delta$ is $g$-invariant, $g$ maps $B$ either onto itself or onto
another ball. Thus if $B \cap a$ is non-empty, $g$ must send $B$ to
itself.

Double $B$ along its boundary to get $S^3$. If $g$ maps $B$ to
itself, its action will extend to an orientation-preserving, periodic,
smooth map $h:S^3 \to S^3$. By the solution to the Smith Conjecture
\cite{smith}, the fixed-point set of $h$ is a single unknotted
circle. Thus the double of $B \cap a$ is the unknot, making $B \cap a$
a single, boundary-parallel arc. We can choose the isotopic arc $b
\subset \bdy B$ so that it intersects the $E_j$ in a minimal way: if
an endpoint of $B \cap a$ lies in some $E_j$, then $b \cap E_j$ is a
radius of that disk. Otherwise, $b$ is disjoint from the
distinguished disks.

\begin{figure}
\psfrag{a}{$a$}
\psfrag{b}{$b$}
\psfrag{bb}{$b'$}
\psfrag{B}{$B$}
\psfrag{BB}{$B'$}
\psfrag{Ej}{$E_j$}
\psfrag{Ejj}{$E_{j'}$}
\begin{center}
\epsfbox{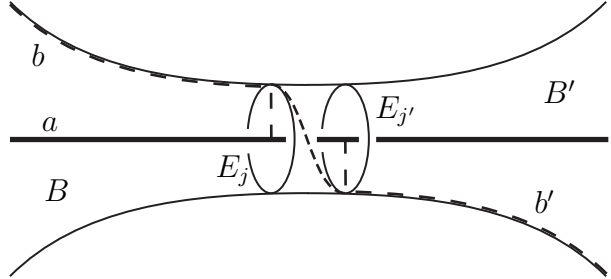}
\end{center}
\caption{When the components of $V \vert \Delta$ are glued together 
along $\Delta$, the fixed-point locus $a$ remains isotopic to the boundary.}
\label{ball-joining}
\end{figure}

To prove that $a$ is boundary-parallel in all of $V$, it remains to
show how the isotopies in the individual balls extend across the disks
in $\Delta$. If some $D_i \in \Delta$ is disjoint from $a$, it
presents no problem. If $D_i$ does intersect $a$, it is preserved by $g$,
and the intersection must be a point or an arc.  If $D_i$ intersects
$a$ in an arc, then that component of $a$ is already boundary-parallel
through $D_i$. If $D_i$ intersects $a$ transversely, in a point, then
the corresponding pair of disks $E_j \subset \bdy B$ and $E_{j'}
\subset \bdy B'$ will each intersect it in a point too. In this case,
the isotopies of $a \cap B$ and $a \cap B'$ to $\bdy V$ can be joined
over $N(D_i)$, as in Figure \ref{ball-joining}. Thus any non-empty
fixed-point locus $a$ must be boundary-parallel.

To prove part $(2)$ of the theorem, consider the quotient of $V \vert
\Delta$ under the action of $g$. This quotient is still a disjoint 
union of balls, which are glued along disks on their boundaries to
reconstruct $W$. Thus $W$ can be viewed as a thickened graph, whose
vertices are in the disjoint balls and whose edges correspond to
gluings. As a thickened graph, $W$ must be a handlebody, and the image
of $a$ under the quotient map is boundary-parallel by the same
argument as above.
\end{proof}

\section{Proof of Theorem \ref{main}}

Lemma \ref{known-case} proves the ``if'' direction of the theorem. To
prove the ``only if'' direction, suppose that a strong involution
$\psi$ of $K$ fixes its unknotting tunnel $\tau$ pointwise. By the
Smith Conjecture solution, the fixed-point locus of $\psi$ is
an unknotted circle. $\tau$ already lies on this axis; call the
remaining arc of the axis $a$.

\begin{notation}
Let $\pi: S^3 \to S^3$ be the quotient map induced by the action of
$\psi$. Label $\Khat = \pi(K)$, $\tauhat = \pi(\tau)$, and $\ahat =
\pi(a)$. Then $\pi$ is a branched covering map of $S^3$ by $S^3$,
branched along the unknot $\tauhat \cup \ahat$.
\end{notation}

Recall the genus-$2$ Heegaard splitting of $S^3$ by $V_1=\overline{N(K\cup 
\tau)}$ and $V_2 = E(K \cup \tau)$. $\psi$ acts as an involution on each 
$V_i$; let $W_i = \pi(V_i)$ be the quotients. Since $W_1$ is a closed
regular neighborhood of the knot $\Khat \cup \tauhat$, it is a solid
torus. By Theorem \ref{smith-hand}, $W_2$ is a handlebody, and since
$T = \bdy W_1 = \bdy W_2$ is a torus, $W_2$ is itself a solid
torus. The result now follows in two steps.

\begin{claim}\label{2br-downstairs}
$\Khat \cup \ahat$ is a two-bridge knot with splitting sphere 
$\bdy N(\tauhat)$.
\end{claim}

\begin{proof}
The involution $\psi$ acts on each $V_i$ separately, and $a$ is the
fixed-point set of $\psi$ in $V_2$.  By Theorem \ref{smith-hand}, it
follows that $\ahat$ is boundary-parallel in $W_2$. Thus $W_2
\setminus N(\ahat) = E(\Khat \cup \tauhat \cup \ahat)$ is a genus-$2$ 
handlebody, and $\tauhat$ is an unknotting tunnel for $\Khat \cup
\ahat$. Furthermore, $\tauhat \cup \ahat$ is the unknot by the solution 
to the Smith Conjecture, and $\Khat \cup \tauhat$ is the unknot
because $W_2 = E(\Khat \cup \tauhat)$ is a solid torus. Thus Corollary
\ref{recognition} tells us that $\bdy N(\tauhat)$ splits $\Khat \cup \ahat$ 
into rational tangles.
\end{proof}

\begin{claim}
$K$ is a two-bridge knot, and $\tau$ is its upper or lower tunnel.
\end{claim}

\begin{proof}
Claim \ref{2br-downstairs} tells us that $\bdy N(\tauhat)$ is a
splitting sphere for $\Khat \cup \ahat$. In particular, $\Khat \cap
E(\tauhat)$ is isotopic to $\bdy N(\tauhat)$ via a disk $\hat{D}$ that
is disjoint from $\ahat$. Since $\hat{D}$ is disjoint from the branch
locus of $\pi$, it lifts to two disjoint disks, $D_1$ and $D_2$, that
realize isotopies of $K_1$ and $K_2$, respectively, to $\bdy N(\tau)$.
Therefore, $(E(\tau), K \cap E(\tau))$ is a rational tangle and $K$ is
a $2$-bridge knot.

It follows that $\tau$ is an upper or lower tunnel for $K$, by the
same argument as in Corollary \ref{recognition}.
\end{proof}

\bibliographystyle{plain} 

\bibliography{involution}

\end{document}